\documentclass{ifacconf}
 
\usepackage{amsmath,amssymb, bm, mathtools} 
\usepackage{graphics} 
\usepackage{color}
\usepackage{float}
\usepackage{comment}
\usepackage{tcolorbox} 
\usepackage{natbib}        
 
\usepackage[section]{algorithm}
\usepackage{algorithmicx}
\usepackage{algpseudocode}
\newcommand{\Input}{\State \textbf{Input: }}  
\newcommand{\Output}{\State \textbf{Output: }}  
\newcommand{\Init}{\State \textbf{Init: }}  

\newcommand{\E}{\mathbb{E}}
\newcommand{\tr}{\text{tr}}

\begin{document}
\begin{frontmatter}
\small \textcopyright 2026 R.P.P.F. Goetz, Y. Dwaraga, N. van de Wouw, T. Oomen, M.M.J. van de Wal, B. Sharif, and H.J. Zwart. This work has been accepted to IFAC World Congress 2026 for publication under a Creative Commons Licence CC-BY-NC-ND.
\vspace{-.4cm}
\title{Optimal Sensor Placement for Output Estimation Using an Artificial Bee Colony Algorithm with Pre-filter \thanksref{footnoteinfo}} 

\thanks[footnoteinfo]{This research is part of the project ThermOpt D$\&$C that  is co-financed by Holland High Tech, top sector High-Tech Systems and Materials, with a PPP innovation subsidy for public-private partnerships for research and development.}

\author[First]{R.P.P.F. Goetz} 
\author[First]{Y. Dwaraga} 
\author[First]{N. van de Wouw}
\author[First]{T. Oomen}
\author[Third]{M.M.J. van de Wal}
\author[Third]{B. Sharif}
\author[First,Second]{H.J. Zwart}

\address[First]{Department of Mechanical Engineering, Eindhoven University of Technology, 5600 MB Eindhoven, The Netherlands (e-mail: r.p.p.f.goetz@tue.nl, y.e.dwaraga@tue.nl, n.v.d.wouw@tue.nl, t.a.e.oomen@tue.nl, h.j.zwart@tue.nl)}
\address[Second]{Department of Applied Mathematics, University of Twente, 7500 AE Enschede, The Netherlands}
\address[Third]{ASML, 5504 DR Veldhoven, The Netherlands, (e-mail: marc.van.de.wal@asml.com, bardia.sharif@asml.com)}

\begin{abstract}  
    Sensor placement for maximizing the estimation performance of the Kalman filter is an NP-hard optimization problem. Furthermore, its feasible set grows combinatorially with the candidate locations and the number of sensors. In this paper, we study this sensor placement problem for a 3D thermoelastic system modelled as a discrete-time linear stochastic model. We use the Novel Binary Artificial Bee Colony (NBABC) algorithm  with a Gramian-based pre-filter to reduce the computational complexity. Our results show the efficiency and the fast convergence of the proposed approach.
\end{abstract}

\begin{keyword}
    Sensor placement, Bio-inspired algorithms and optimization-based control, Estimation and filtering, Optimization-based estimation and control, High-performance motion control systems, Linear systems, Kalman filtering.
\end{keyword}

\end{frontmatter}

\section{INTRODUCTION}

In control systems, estimating unmeasurable quantities is often necessary, and the resulting estimation quality is highly dependent on the placement of sensors. For high-dimensional systems such as structural, see \cite{guptaOptimizationCriteriaOptimal2010}, or thermal applications, solving the sensor placement problem is challenging as many locations are possible. In this work, we focus on sensor placement problems occurring in high-precision lithography machines used for chip manufacturing. Such machines produce wafers on which electronic patterns are printed using Ultraviolet (UV) light. Since a chip is composed of several layers, the UV light must expose at exactly the desired location. However, heat-induced deformation of the stage moving the wafer hinders such accurate positioning for exposure. Therefore, temperature sensors must be positioned optimally for accurate estimation of this unmeasurable deformation. This accurate estimation is only required for a finite-time interval corresponding to the exposure duration. The estimate can later be used for correcting the positioning. 

The sensor placement objective can be formulated in many ways, see for an overview \cite{guptaOptimizationCriteriaOptimal2010}, depending on the objective pursued. As we focus on discrete-time linear-time-invariant (LTI) stochastic models, the optimal estimator is the Kalman filter. The sensor placement objective is then formulated as the maximization of the estimation performance of the Kalman filter over a given time interval, i.e., the exposure duration. More precisely, the average trace of the finite-horizon estimation error covariance matrix is minimized, as in \cite{goetzOptimalCoDesignSensor2025}. This problem is challenging to solve for several reasons:
\begin{itemize}
    \item The problem is NP-hard, see \cite{zhang_sensor_2017};
    \item The cost function lacks modularity properties, see \cite{jawaidSubmodularityGreedyAlgorithms2015};
    \item The feasible set grows combinatorially.
\end{itemize}

Several approaches have been proposed to solve this sensor placement problem. The first approach is to consider surrogate objectives that exhibit modularity characteristics. Such modularity property enables the use of a greedy algorithm as it converges with closeness guarantees to the global optimum, see \cite{zhang_sensor_2017} or \cite{jawaidSubmodularityGreedyAlgorithms2015}. While this strategy can be computationally cheap, the modularity property only guarantees the greedy algorithm to find a sensor configuration close to the optimum of the surrogate objective and not to the one minimizing the average trace of the covariance matrix. A second approach is to use convexification techniques that enable the use of efficient convex optimization solvers, see \cite{darivandiAlgorithmLQOptimal2013} for the dual problem of actuator placement. However, convexification involves relaxation making the obtained solution suboptimal. A last approach is to use metaheuristics, such as genetic algorithms or swarm optimizations, see \cite{wrobelOptimalSensorPlacement2023}. While these methods may find the global optimum of non-modular NP-hard problems, they remain computationally expensive and offer no guarantees of global optimality.

In this paper, we focus a metaheuristic approach and develop a fast algorithm that can be used for solving high-dimensional sensor placement problems. In particular, we will exploit the recently developed Novel Binary Artifical Bee Colony (NBABC) which is a probabilistic algorithm for solving high-dimensional binary optimization problems, see \cite{santanaNovelBinaryArtificial2019}. To reduce the combinatorial growth of the feasible space and to accelerate the convergence of this algorithm, we introduce a pre-filter screening the candidate sensor locations by keeping the ones satisfying an observability condition. 

The main contribution of this work is the adaptation of the NBABC to solve high-dimensional sensor placement problems and to show its improved convergence when a Gramian-based pre-filter is used.

The article is structured as follows. In Section \ref{sec_ProbDescr}, the sensor placement problem is formulated for discrete-time LTI stochastic models. The strategy used for solving this problem is presented in Section \ref{sec_NBABC}. This strategy is composed of a Gramian-based filter and a metaheuristic optimizer. Finally, this approach is applied in Section \ref{sec_res3D} for finding the best sensor configuration to estimate the heat-induced deformation in a 3D thermoelastic model representing a wafer stage in lithography machines.

\textit{Notations:} A symmetric matrix $A = A^\top \in\mathbb{R}^{n\times n}$ is called positive definite if $x^\top A x > 0$ for all $x\in\mathbb{R}^n\setminus{\{0\}}$, and is written as $A\succ 0$. Similarly, the matrix $A$ is positive semi-definite if $x^\top A x \geq 0$ for all $x\in\mathbb{R}^n$, and is written as $A\succeq 0$. A Gaussian random variable $X$ with mean $\mu$ and variance matrix $\sigma$ is written as $X \sim \mathcal{N}(\mu,\sigma)$. The expected value of the random variable $X$ is denoted by $\mathbb{E}[X]$. The signal $\hat{X}_{s|t}$ denotes an estimate of $X_s$ at time $s$ using measurements until time $t$.
\section{PROBLEM DESCRIPTION}\label{sec_ProbDescr}
\begin{figure}[bp]
    \centering
    \includegraphics[width = 0.85\columnwidth]{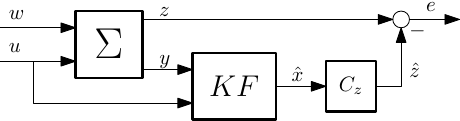}
    \caption{Estimation procedure of $\hat{z}$ using the Kalman filter (KF) and the linear map $C_z$.}
    \label{fig_Est_KF}
\end{figure}

For discrete-time LTI systems, the sensor placement problem is formulated as the maximization of the estimation performance. The estimation procedure is illustrated in Figure \ref{fig_Est_KF}. The physical dynamics under study is described by the model $\Sigma$ that is affected by a known input $u$ and unknown disturbance $w$. The unmeasurable output $z$ must be estimated. This is achieved by the Kalman filter (KF) that uses measurement $y$ and input $u$ to estimate the state $x$ of $\Sigma$. Then an output map $C_z$ is applied on this estimate $\hat{x}$ to obtain the desired output estimate $\hat{z}$. The estimation error $e = z - \hat{z}$ should be minimized by selecting the best sensor layout.

The model $\Sigma$ belongs to the class of stable linear-time-invariant stochastic discrete-time models:
\begin{equation}\label{eq_sys_stocha}
        \Sigma :=
    \begin{bmatrix}
        x_{k+1}\\
        y_k\\
        z_k
    \end{bmatrix}=
    \begin{bmatrix}
        A & B_u & B_w & 0\\
        C & 0 & 0 & I\\
        C_z & 0 & 0 & 0
    \end{bmatrix}
    \begin{bmatrix}
        x_k\\
        u_k\\
        w_k\\
        v_k
    \end{bmatrix}
\end{equation}
with the discrete time instant $k$. The sensor configuration is represented by the matrix $C$ which is restricted to the set $\Omega_C$. The unknown disturbance $w$ and the measurement noise $v$ in \eqref{eq_sys_stocha} are zero-mean Gaussian white noises and the initial state $x_0$ is a Gaussian random variable:
\begin{equation}\label{eq_as_wvx0}
        w_k\sim\mathcal{N}(0,Q),\quad v_k\sim\mathcal{N}(0,R),\quad x_0\sim\mathcal{N}(\hat{x}_{0|0},P_{0|0})
\end{equation}
with known symmetric covariance matrices $Q\succeq 0$, $R\succ 0$, $P_{0|0}\succeq 0$. These signals are considered mutually independent.

The sensor placement problem entails maximizing the transient (or finite-horizon) output estimation performance of the Kalman filter over the time interval $[0,k_f]$ corresponding to the exposure duration in the considered application:
\begin{equation}\label{eq_FH_SP}
    \min_{C\in\Omega_C} \frac{1}{k_f}\sum_{k=1}^{k_f}\tr[P_{z_{k|k}}],
\end{equation}
where the output estimation error covariance is given by 
\begin{align}
    P_{z_{k|k}} &= \E \left[(z_k-\hat{z}_{k|k})(z_k-\hat{z}_{k|k})^\top\right]\\
    & = C_z\E \left[(x_k-\hat{x}_{k|k})(x_k-\hat{x}_{k|k})^\top\right] C_z^\top\\
    &= C_z P_{k|k}C_z^\top
\end{align}
with $P_{k|k}$ being the a posteriori state covariance matrix, see \cite{gelbAppliedOptimalEstimation1974} Chapter 4. As $P_{k|k}$ depends on $C$, the sensor layout directly impacts the desired output estimation performance.

The sensor placement problem \eqref{eq_FH_SP} for the system in Figure \ref{fig_Est_KF} is challenging for the three reasons mentioned in the introduction. Besides the NP-hardness and non-modularity, the feasible set grows combinatorially as its size equals the total number of possible sensor configurations given by
\begin{equation}\label{eq_n_conf}
    \frac{n_{cand}!}{n_y!(n_{cand} - n_y)!},
\end{equation}
where $n_y$ is the number of sensors to be positioned among $n_{cand}$ candidate locations. These challenges motivate the design of a computationally efficient algorithm capable of solving the sensor placement problem with high accuracy and reasonable computational cost.

We focus on binary $C$ matrices where each row corresponds to a sensor. When the state $x\in \mathbb{R}^n$  represent physical variables at $n$ discrete spatial locations (e.g., temperature at the nodes of a finite-element mesh for the thermoelastic problem considered in this paper), the entry $C_{ij}=1$ indicates that the sensor $i$ measures the state at the $j$th node whereas $C_{ij}=0$ indicates the absence of that sensor. Note that only point sensors are considered meaning that sensors only measure at one node. With this structure, the problem becomes a combinatorial optimization problem which justifies the study of binary variants of metaheuristic algorithms.

\section{THE NOVEL BINARY ARTIFICIAL BEE COLONY ALGORITHM WITH PRE-FILTERING}\label{sec_NBABC}
Because of the complexity of the optimization problem \eqref{eq_FH_SP}, we focus on a probabilistic search algorithm and especially on the Artificial Bee Colony (ABC) algorithm initially proposed by \cite{karabogaIdeaBasedHoney2005}, which demonstrated superior performance over other metaheuristics for a large set of numerical test functions, see \cite{karabogaComparativeStudyArtificial2009}. A binary version of the algorithm is selected and adapted for the sensor placement problem. Moreover, to address the combinatorial explosion of the feasible set, we introduce a Gramian-based condition to pre-select sensor locations exhibiting sufficient observability properties.

\subsection{The NBABC algorithm for sensor placement}
The Novel Binary Artificial Bee Colony (NBABC) algorithm is a binary variant of the original ABC. Many variants were developed to improve its performance, see the latest survey article by \cite{ibrahimArtificialBeeColony2025} and the article by \cite{akaySurveyArtificialBee2021} focusing on binary, integer and mixed-integer problems. Among these variants, we adopt the NBABC algorithm proposed by \cite{santanaNovelBinaryArtificial2019} for several reasons. First, this variant is binary as required by our optimization problem \eqref{eq_FH_SP}. Second, its crossover mechanism for generating new candidate solutions uses logic gates which are computationally inexpensive. Third, the authors \cite{santanaNovelBinaryArtificial2019} demonstrated that the NBABC outperforms nine other binary algorithms on high-dimensional and NP-hard problems.

In this section, we detail the NBABC algorithm and adapt it for solving the sensor placement problem. 

\subsubsection{The NBABC as a bio-inspired algorithm}\hfill \\
Similar to other ABC variants, the NBABC is a multi-agent probabilistic search algorithm bio-inspired from the behavior of honey bees. In a honey bee colony, there exist three types of bees:
\begin{itemize}
    \item The employed bee (EB), associated to a food source, explores a neighboring source and dances in the hive based on the quality of this new source found,
    \item The onlooker bee (OB) probabilistically decides whether to follow an EB based on the dance and explores the neighborhood of the associated food source,
    \item The scout bee (SB) randomly explores to find a new food source once a known food source is abandoned.
\end{itemize}
The general structure of the NBABC is given in Algorithm \ref{algo_Gen_NBABC}. Usually, an EB is associated with one food source and the total number of EB defines the population size which also equals the number of OB. Both the EB and OB aim to improve the exploited food sources. If a food source is not improved after some attempts, the associated EB abandons the source and becomes a SB (Line \ref{alg_l_SB} in Algorithm \ref{algo_Gen_NBABC}). After randomly exploring around the hive, the SB locates a new food source and resumes to an EB. The mechanism of both the EB and OB to find a new food source (Lines \ref{alg_l_EB_Find} and \ref{alg_l_OB_Find} in Algorithm \ref{algo_Gen_NBABC}) relies on the food sources present in the population. This crossover mechanism is detailed later in this section.

\begin{algorithm}[btp]
    \caption{General structure of the NBABC Algorithm}
    \label{algo_Gen_NBABC}
    \begin{algorithmic}[1]
        \Input Population size 
        \Init Initialize randomly the food sources of the colony\label{alg_l_init}
        \While {Termination condition is false}\label{alg_l_termin}
        \For {Each EB}
        \State Finds a new food source in the neighborhood\label{alg_l_EB_Find}
        \State Evaluates the quality of the new source
        \State Keeps the new or old source with better quality\label{alg_l_EB_select}
        \EndFor
        \State Compute the probability that EBs attract OBs\label{alg_l_proba}
        \For {Each OB}
        \State Follows an EB based on the probability
        \State Finds a new food source in the neighborhood\label{alg_l_OB_Find}
        \State Evaluates the quality of the new source
        \State Keeps the new or old source with better quality\label{alg_l_OB_select}
        \EndFor
        \If {A food source is abandoned}\label{alg_l_SB}
        \State The associated EB becomes SB
        \State Find randomly a new food source
        \State The SB resumes to EB
        \EndIf
        \State Update the best food source found so far
        \EndWhile
        \Output Best food source found so far
    \end{algorithmic}
\end{algorithm}

\subsubsection{Adaptation for sensor placement}\hfill \\
In terms of the sensor placement problem, a food source corresponds to a feasible solution of the optimization problem, which is a sensor configuration $C\in\Omega_C$. The binary matrix $C\in\{0,1\}^{n_y\times n}$ can be equivalently expressed as a binary vector $c\in\{0,1\}^{n_{cand}}$ where the $j$th entry equals $1$ if a sensor is placed at the $j$th candidate node and $0$ otherwise. Because the number of sensors is fixed to $n_y$, $c$ contains exactly $n_y$ ones. The quality of a sensor configuration $c_i$, $i\in\{1,2,..., n_{pop}\}$, is called the fitness. Its standard definition is used:
\begin{equation}\label{eq_ABC_fitness}
    fit(c_i) = \frac{1}{1+f(c_i)},
\end{equation}
where $f=\frac{1}{k_f}\sum_{k=1}^{k_f}\tr[P_{z_{k|k}}]$ is the cost function in \eqref{eq_FH_SP}. Every time the cost is evaluated, a counter is incremented. The algorithm terminates (Line \ref{alg_l_termin} in Algorithm \ref{algo_Gen_NBABC}) when this counter exceeds a user-defined threshold $N_f$. The probability of an OB to follow the EB associated to the solution $i$ (Line \ref{alg_l_proba} in Algorithm \ref{algo_Gen_NBABC}) is given by
\begin{equation}\label{eq_ABC_pi}
    p_i = \frac{fit(c_i)}{\sum_{j=1}^{n_{pop}}fit(c_j)}.
\end{equation}

The EB and OB phases (Lines \ref{alg_l_EB_Find} and \ref{alg_l_OB_Find} in Algorithm \ref{algo_Gen_NBABC}) are adapted for sensor placement. In the EB phase, the current solution $c_i$ is the configuration associated to the EB. In the OB phase, the bee selects a solution $c_i$ with probability $p_i$ in \eqref{eq_ABC_pi}. A new sensor configuration is generated by a crossover mechanism that preserves the fixed number of sensors $n_y$. First, a partner solution $c_j$ is randomly chosen from the set of solutions of the population. The positions at which $c_i$ and $c_j$ differ are identified. From these differing positions, a user-defined number of flips $N_{0\leftrightarrow 1}$ is chosen: $N_{0\leftrightarrow 1}$ entries where $c_i = 1$ are randomly selected and flipped to $0$, and an equal number where $c_i=0$ are randomly selected and flipped to $1$. Applying these flips produces a new candidate sensor configuration. A trial counter tracks the number of unsuccessful improvement attempts in both phases. If the newly generated solution has a better fitness than the current one, the EB exploits this improved solution and the trial counter is reset to 0. Otherwise, the current solution is retained, and its trial is incremented (Lines \ref{alg_l_EB_select} and \ref{alg_l_OB_select} in Algorithm \ref{algo_Gen_NBABC}). When the counter reaches the limit $N_{\mathrm{trial}}$, the solution is abandoned and the SB phase is activated (Line \ref{alg_l_SB} in Algorithm \ref{algo_Gen_NBABC}). 

In the algorithm, additional parameters have been introduced: the number of flips $N_{0\leftrightarrow 1}$, the maximal improvement attempts $N_{\mathrm{trial}}$, the maximum function evaluation $N_f$. Naturally, the overall performance of the NBABC depends on proper tuning of these parameters.



\subsection{Gramian-based pre-filter}
To limit the combinatorial growth of the feasible solution space, we pre-filter the possible sensor locations with an observability condition before applying the NBABC. It is expected that configurations that produce good estimation performance also give good observability of the system\footnote{Note that the configurations with the best observability properties do not necessarily relate to the optimal sensor location of \eqref{eq_FH_SP}.}. Therefore, a threshold can be used to remove locations with poor observability. The observability measure is inspired by the formulation proposed in \cite{merks_towards_2019}, Section 3.2, where it is directly maximized for sensor placement. As explained below, a sensor configuration should provide good observability trade-off between the contribution of the performance output $z$ and that of the disturbance $w$, leading to the Gramian-based objective:
\begin{equation}\label{eq_obs_criterion_notnormed}
        \tilde{h}_\mathcal{O}(C)=\alpha_1 \tr(C_z \mathcal{W} C_z^\top) + \alpha_2 \tr(Q^{1/2}B_w^\top \mathcal{W} B_w Q^{1/2})
\end{equation}
with weights $\alpha_1$ and $\alpha_2$. The first term in \eqref{eq_obs_criterion_notnormed} represents the observability associated to the output $z$. In fact, the state observability Gramian, 
\begin{equation}\label{eq_state_obs}
    \mathcal{W} = \sum_{k=0}^{\infty} A^{k^\top} C^\top CA^k,
\end{equation}
is projected onto the space spanned by $C_z$. The second term in \eqref{eq_obs_criterion_notnormed} gives the observability regarding an impulse in $w$ (with $x_0 = 0$). The energy in the measurement $y$ is computed using the associated system response. By scaling $w = Q^{1/2} \tilde{w}$ and having an impulse induced by $\tilde{w}$, we have
\begin{align}
    \sum_{k=0}^{\infty} y_k^\top y_k &= \sum_{k=0}^{\infty}\left[\left(\sum_{j=0}^{k-1}\tilde{w}_j^\top Q^{1/2^\top}B_w^\top A^{{k-j-1}^\top} C^\top\right)\right.\nonumber\\
    &\quad\left. \left(\sum_{i=0}^{k-1} C A^{k-i-1} B_w Q^{1/2} \tilde{w}_i\right)\right]\nonumber\\
    &=\sum_{k=1}^{\infty}Q^{1/2^\top}B_w^\top A^{{k-1}^\top} C^\top C  A^{k-1}B_w Q^{1/2}\nonumber\\
    &=  Q^{1/2^\top}B_w^\top \mathcal{W} B_w Q^{1/2},
\end{align}
which relates to the second term in \eqref{eq_obs_criterion_notnormed}. This energy directly indicates how sensitive the measurements are to an impulse in $\tilde{w}$. Because the two terms in \eqref{eq_obs_criterion_notnormed} could actually have different order of magnitude, each term is scaled:
\begin{equation}\label{eq_obs_criterion_normed}
    \begin{aligned}
        h_{\mathcal{O}}(C) =\, &\alpha \frac{\tr(C_z \mathcal{W} C_z^\top)}{\max_C \tr(C_z \mathcal{W} C_z^\top)} \\
        &+ (1-\alpha) \frac{\tr(Q^{1/2}B_w^\top \mathcal{W} B_w Q^{1/2})}{\max_C\tr(Q^{1/2}B_w^\top \mathcal{W} B_w Q^{1/2})}
    \end{aligned}
\end{equation}
with $0\leq\alpha \leq 1$ to ensure $0\leq h_{\mathcal{O}}(C)\leq1$.

Due to the modularity property of the Gramian $\mathcal{W}$ with respect to $C$, see \cite{summersSubmodularityControllabilityComplex2016}, it is sufficient to check the observability metric \eqref{eq_obs_criterion_normed} per individual sensor location instead of the (many) complete configurations. Then, the potential sensor locations that do not satisfy the condition $h_{\mathcal{O}}>\bar{h}_{\mathcal{O}}$ are removed, where $\bar{h}_{\mathcal{O}}$ is a threshold value defined by the user. This reduces the number $n_{cand}$ of potential sensor locations and, with it, the number of possible sensor configurations in \eqref{eq_n_conf}, but it still carries the risk of removing optimal locations. 
\section{CASE STUDY RESULTS}\label{sec_res3D}
\begin{figure}[bp]
    \centering
    \includegraphics[width = \columnwidth]{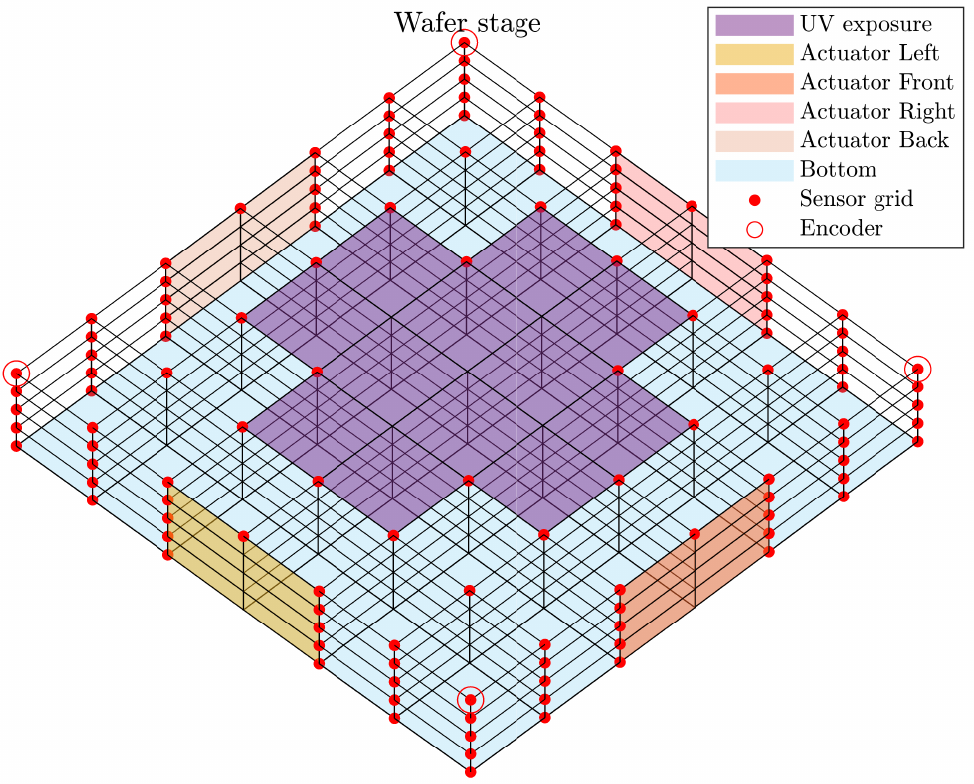}
    \caption{Case study: Discretized wafer stage}
    \label{fig_setup}
\end{figure}

As a case study, we consider a stochastic model representing the thermoelastic dynamics of a wafer stage in lithography machines, see Figure \ref{fig_setup}. The dynamics is discretized in space and time. The state $x$ represents the temperature field, $z$ is the heat-induced displacement at the position encoders, $y$ is the temperature measurements affected by noise $v$, the known heat input is $u$ and the unknown heat disturbance is $w$. Without loss of generality, we take $u=0$. Six heat disturbances are considered: heat coming from the UV exposure, from the motion table below the wafer stage and from the four motion actuators. For further details and parametric settings, see \cite{goetzOptimalCoDesignSensor2025}. The model, the algorithm and the analyses were implemented in MATLAB (The MathWorks, Inc.).

For the analyses, we aim to place two relatively precise sensors, $n_y=2$ among $n_{cand} = 124$ candidate locations, with noise $R = 10^{-6}\, \mathrm{K^2}$. Moreover, the uncertainty of the initial temperature field is captured by $P_{0|0}$. We suppose this matrix to be diagonal with perfect knowledge at the actuator nodes making the variances zero. For the other nodes, we suppose to have no perfect knowledge, i.e., the variances equal $0.001 \mathrm{K^2}$. About the disturbance covariance, we fix it to $Q=\mathrm{diag}(([7,1,3,2.5,2,1.5]/\mathrm{a_{Area}})^2)\,(\mathrm{W/m}^2)^2$, where $\mathrm{a_{Area}}$ contains the areas in $\mathrm{m^2}$ of each disturbance ordered as follows: the UV exposure, the bottom table, and the left, front, right and back actuators.

To apply the NBABC, additional parameters must be specified. We assumed a population size of 10 bees. The $N_{\mathrm{trial}}$ is set to $10$. During the crossover step in the EB and OB phases, only one element is allowed to be flipped, i.e., $N_{0\leftrightarrow 1} = 1$. The limit of function evaluations is $N_f = 500$.

\subsection{Convergence with and without pre-filtering}\label{ssec_Res_withVSwithout}
\begin{figure}[btp]
    \centering
    \includegraphics[width = \columnwidth]{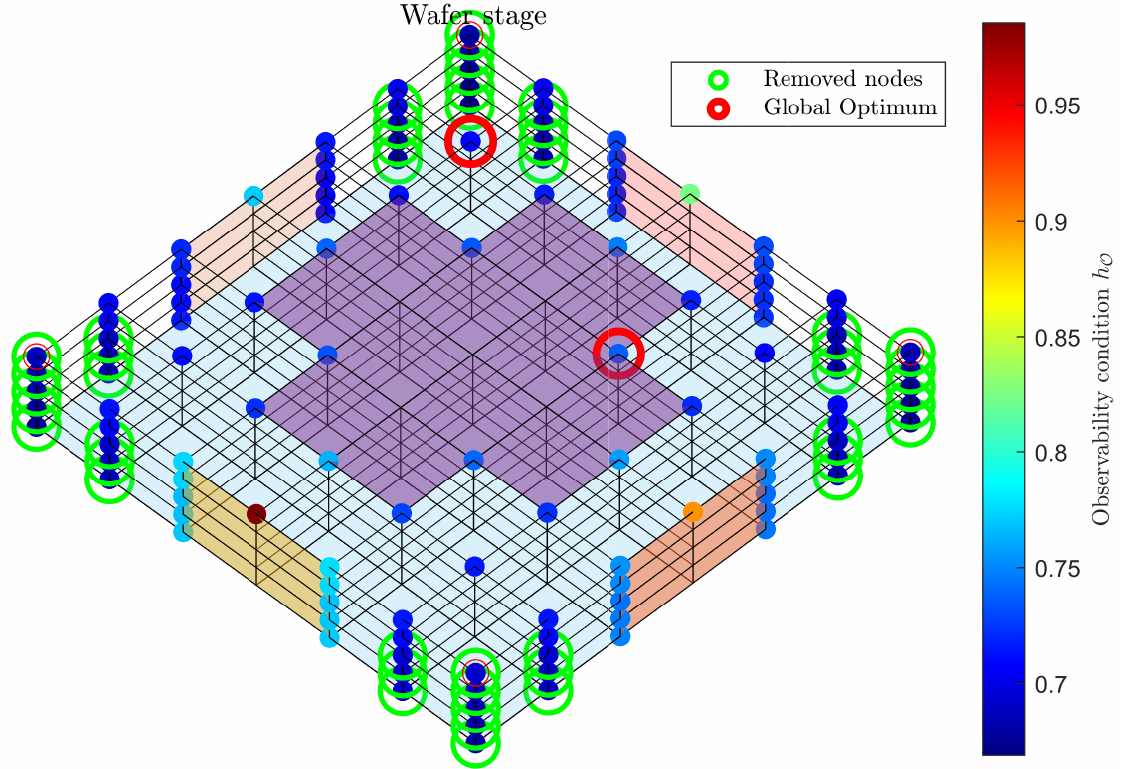}
    \caption{Nodes removed by the filter with $\bar{h}_{\mathcal{O}} = 0.7$}
    \label{fig_Rem_Nodes}
\end{figure}

\begin{figure}[btp]
    \centering
    \includegraphics[width = 0.969\columnwidth]{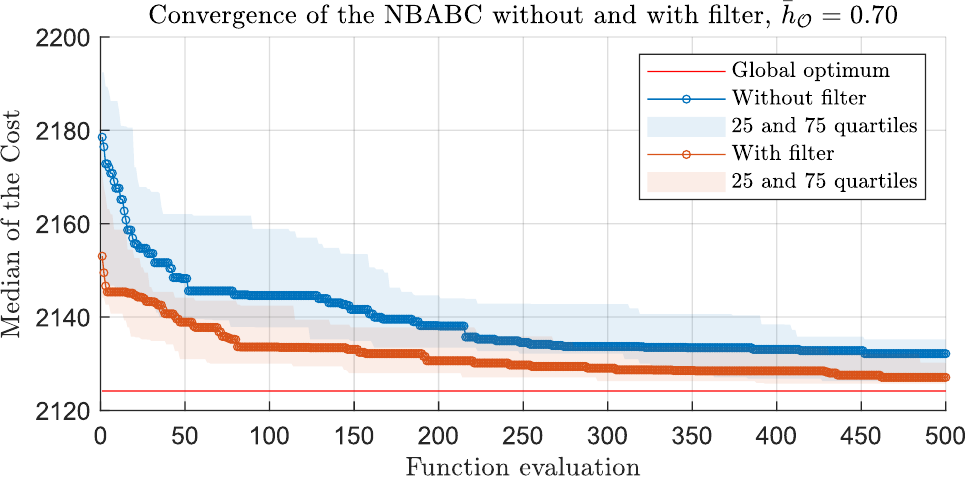}
    \caption{Median of the cost per function evaluation over 50 runs without and with filter using $\bar{h}_{\mathcal{O}} = 0.7$}
    \label{fig_Conv}
\end{figure}

In this section, we compare the convergence of the NBABC with and without pre-selection of the nodes. We set equal importance to both terms in \eqref{eq_obs_criterion_normed}, $\alpha = 0.5$, and use a threshold $\bar{h}_{\mathcal{O}} = 0.7$ to remove 41 out of the 124 initial candidate locations, see Figure \ref{fig_Rem_Nodes}. Because NBABC is probabilistic, the convergence of the algorithm is assessed using the median cost over 50 runs of the NBABC. Figure \ref{fig_Conv} shows this convergence with and without filtering as well as the global optimum obtained by exhaustive search. 

We first observe that the initial cost is lower with filter than without, meaning that the removed nodes are indeed part of high-cost configurations. Next, we see in Figure \ref{fig_Rem_Nodes} that the optimal sensor locations are not removed by the filter with the selected threshold. In Figure \ref{fig_Conv}, the median cost converges faster with filter than without. Moreover, out of the 50 runs, the approach with pre-filter found the global optimum 11 times whereas the NBABC without filter found the global optimum only 7 times.

\subsection{Effect of a high threshold $\bar{h}_\mathcal{O}$ for the pre-filter}\label{ssec_Res_highThresh}
\begin{figure}[btp]
    \centering
    \includegraphics[width = \columnwidth]{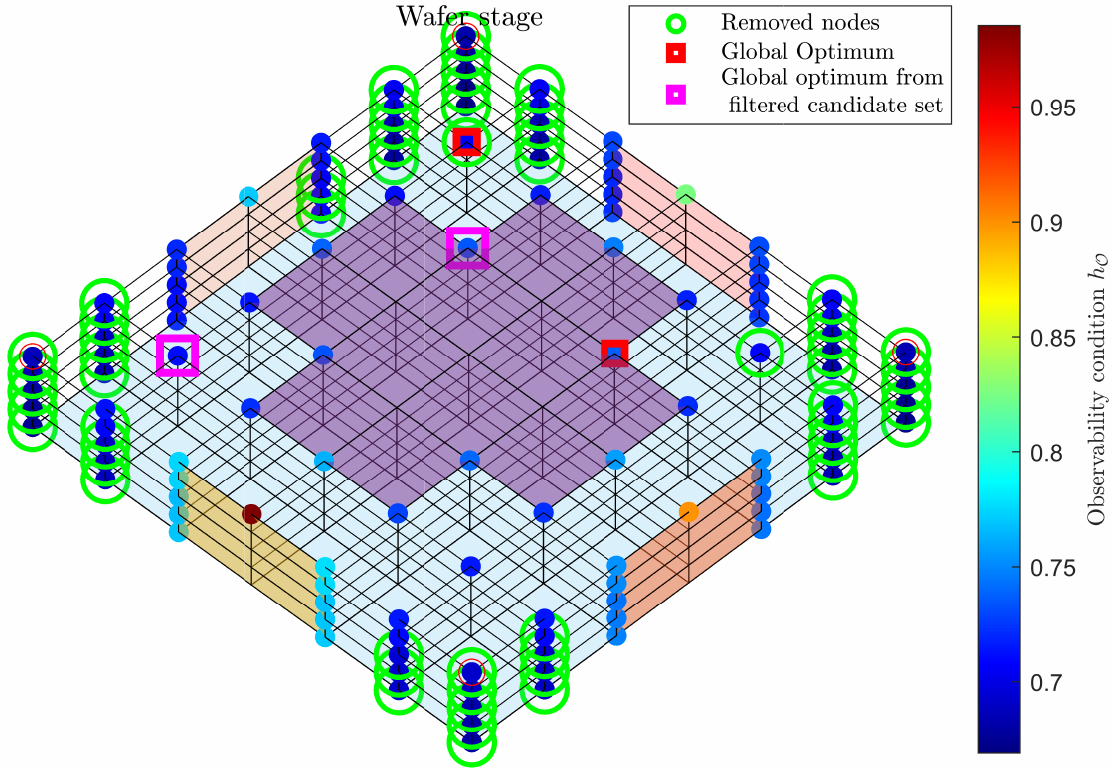}
    \caption{Nodes removed by the filter with $\bar{h}_{\mathcal{O}} = 0.71$}
    \label{fig_Conf_highThresh}
\end{figure}

\begin{figure}[btp]
    \centering
    \includegraphics[width = 0.97\columnwidth]{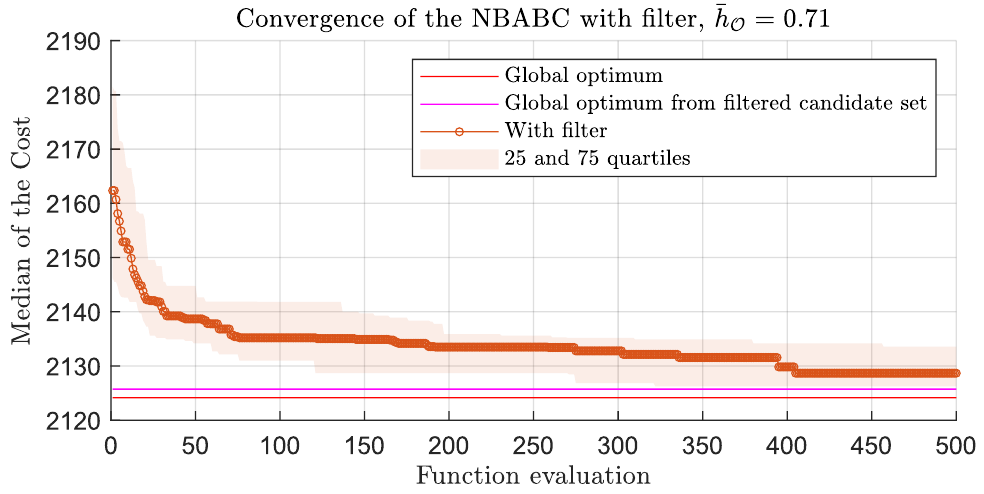}
    \caption{Median of the cost per function evaluation over 50 runs of the NBABC with filter using $\bar{h}_{\mathcal{O}} = 0.71$}
    \label{fig_Conv_highThresh}
\end{figure}

Suppose now that the filter's threshold $\bar{h}_{\mathcal{O}}$ is slightly increased from $0.7$ to $0.71$. In that case, 61 nodes are removed including one of the global optimal locations, see Figure \ref{fig_Conf_highThresh}. In Figure \ref{fig_Conv_highThresh}, the optimum out of the filtered set of feasible configurations is also shown in magenta. We like to stress that even by removing optimal nodes with the filter, the NBABC still finds the best configuration possible from the remaining possible locations. Out of the 50 runs, the optimum was found 10 times.


\subsection{Discussion}
The proposed approach of pre-selecting the candidate nodes via a Gramian-based condition before applying the NBABC shows improved convergence compared to the case without pre-selection. The results show that this approach accelerates the convergence of the NBABC, see Figure \ref{fig_Conv}. So, less function evaluations are needed, making this approach computationally more attractive for solving high-dimensional sensor placement problems. 

A challenge of this approach resides in the parameter tuning. Within the NBABC and the Gramian-based filter, many parameters must be defined, whose values impact the performance. For the ABC, we refer the reader to the article by \cite{akayParameterTuningArtificial2009} for the effect of the different parameters. However, no studies have been performed specifically on the tuning of the NBABC variant. Regarding the Gramian-based pre-selection, a poorly chosen threshold can remove optimal locations, see Figure \ref{fig_Conf_highThresh}. Even if the NBABC still converges to the optimum in the remaining candidate set, this found optimum gives worse estimation performance than the global optimum of \eqref{eq_FH_SP}. One way to avoid discarding optimal locations is to replace the threshold value $\bar{h}_{\mathcal{O}}$ by a specified percentage of candidates to remove. In Section \ref{ssec_Res_withVSwithout}, around $30\%$ of the candidate locations were removed, whereas in Section \ref{ssec_Res_highThresh}, this percentage rises to around $50\%$. We suggest keeping this percentage below $30\%$ to limit the risk of removing optimal locations.

Although maximizing observability does not guarantee good estimation performance as mentioned in \cite{merks_towards_2019}, Chapter 3, the optimal configuration should still induce sufficient observability of the system, motivating the use of this function $h_\mathcal{O}$. The weight $\alpha$ reflects the user's preference: a value close to $1$ emphasizes observability of the initial condition, while a value near $0$ emphasizes observability through the inputs. Other filter functions could be used. One can think for example of computing the estimation performance, i.e.\ the cost function of \eqref{eq_FH_SP}, per node location (as if one sensor is positioned) and select the ones leading to sufficient performance.
\section{CONCLUSIONS AND FUTURE WORK}

In this work, we presented an approach for optimal sensor placement for output estimation that combines a pre-filtering step with the NBABC solving faster the sensor placement problem known to be NP-hard. The sensor placement problem was formulated as the minimization of the average trace of the finite-horizon covariance matrix obtained from the Kalman filter. We studied this problem for stable LTI stochastic discrete-time models. In view of the challenges present in lithography applications, we illustrated the approach for a 3D thermoelastic model representing a wafer stage. Our results show a significant acceleration of the convergence of the NBABC when combined with the proposed pre-filter for removing candidate sensor locations. This strategy makes the sensor placement problem more practical for high-dimensional systems, frequent in industrial applications, by decreasing the computational complexity.

Future work will focus on further improving computational efficiency, in the light of high-dimensional industrial applications. To achieve this, other strategies should be explored. The first idea is to propose a smart initialization of the NBABC population based on engineering insight. This could allow the population to start relatively close to the optimum. Another strategy is to reduce the computational expense of evaluating the objective function. Currently, a single run of the NBABC takes about $20.79\,\mathrm{s}$ from which $98\%$ is spent on cost function evaluations. Rather than evaluating the full objective function, we could instead use a cheaper approximation of this cost.


\bibliography{biblio} 
  
\end{document}